\documentclass[12pt,thmsa]{article}
\usepackage{amsmath, latexsym, amsfonts, amssymb, amsthm, amscd}

\newtheorem{theorem}{Theorem}
\newtheorem{corollary}{Corollary}
\newtheorem{proposition}{Proposition}
\newtheorem{lemma}{Lemma}

\newcommand{\p}{\Bbb{P}}
\newcommand{\px}{\Bbb{P}_x}
\newcommand{\e}{\Bbb{E}}
\newcommand{\ex}{\Bbb{E}_x}

\newcommand{\ind}{\mbox{\rm 1\hspace{-0.04in}I}}

\newcommand{\R}{\mbox{\rm I\hspace{-0.02in}R}}

\newcommand{\ud}{\mathrm{d}}
\newcommand{\ed}{\stackrel{(d)}{=}}

\def\QED{\hfill\vrule height 1.5ex width 1.4ex depth -.1ex \vskip20pt}

\begin{document}

\title{On continuous state branching processes: conditioning and self-similarity.}

\maketitle

\begin{center}
 {\large  A.E. Kyprianou\footnote{\noindent$^{,2}$ Department of Mathematical Science, University of Bath. {\sc Bath, BA2 7AY. United Kingdom.}
$^{1}$E-mail: a.kyprianou@bath.ac.uk, \,$^2$E-mail: jcpm20@bath.ac.uk.\\
$^*$ Corresponding author.}$^{,*}$ and J.C.
Pardo$^2$}
\end{center}
\vspace{0.2in}

\begin{abstract}
In this paper, for $\alpha\in(1,2]$, we show that the $\alpha$-stable continuous-state branching processes and the associated process conditioned never to become extinct are positive self-similar Markov processes. Understanding the interaction of the Lamperti transformation for 
continuous state branching processes and the Lamperti transformation for positive self-similar Markov processes permits access
to a number of explicit results concerning the paths of stable-continuous state branching processes and its conditioned version.

\bigskip


\end{abstract}

\vspace{0.5cm}

\section{Introduction.}
The purpose of this paper is to study self-similarity properties of continuous state branching processes and their conditioned versions when driven by a spectrally positive  $\alpha$-stable process with $\alpha\in(1,2]$. In doing so, a number of results of an explicit nature will fall out of our analysis.
We begin by recalling a number of fundamental facts concerning the stochastic processes that that are of primary interest in this article.

\subsection{Spectrally positive L\'evy processes} 
Let $(\p_x,\, x\in \R)$ be a family of probability measures on the
Skorokhod space, denoted $\mathcal{D}$, such that for each $x\in\R$, the
canonical process $X$ is a L\'evy process with no negative jumps
issued from $x$. Set $\p:=\p_0$, so $\p_x$ is the law of $X+x$
under $\p$. The Laplace exponent $\psi:[0,\infty)\to
(-\infty,\infty)$ of $X$ is specified by
$\e(e^{-\lambda X_t})=e^{t\psi(\lambda)}$, for $ \lambda\in\R$,
and can be expressed in the form
\begin{equation}\label{lk}
\psi(\lambda)=a\lambda+\beta
\lambda^2+\int_{(0,\infty)}\big(e^{-\lambda x}-1+\lambda x\ind_{\{x<1\}}\big)\Pi(\ud x),
\end{equation}
where $a\in \R$, $\beta\geq 0$ and $\Pi$ is a
$\sigma$-finite measure such that
\[
\int_{(0,\infty)}\big(1\land x^2\big)\Pi(\ud x)<\infty.
\]

Henceforth, we shall assume that $(X,\p)$ is not a subordinator (recall that  a subordinator is a L\'evy process with increasing sample paths).
In that case, it is known that the Laplace exponent $\psi$ is strictly convex 
 and tends to $\infty$ as $\lambda\uparrow\infty$. In this case, we define for $q\geq 0$
\[
\Phi(q)=\inf\big\{\lambda\geq 0: \psi(\lambda)>q\big\}
\]
the right inverse of $\psi$ and then $\Phi(0)$ is the largest root of the 
equation $\psi(\lambda)=0$. Theorem VII.1 in \cite{Be} implies that condition $\Phi(0)>0$ holds
if and only if the process drifts to $\infty$. Moreover, almost surely, the paths of $X$
 drift to $\infty$, oscillate or drift to $-\infty$ accordingly 
as $\psi'(0+)<0$, $\psi'(0+)=0$ or $\psi'(0+)>0$.

\subsection{Conditioning to stay positive}

In this article, we also deal with L\'evy processes conditioned to
stay positive. The following commentary is taken from  Chaumont and Doney \cite{CD}
and Chapter VII of Bertoin \cite{Be}. 
In the current context, the L\'evy process, $X$,
conditioned to stay positive is the strong Markov process whose law is given by
\begin{equation}
 \p^\uparrow_x  (X_{t}\in dy) = \lim_{q\downarrow 0}\px(X_t \in dy, \, t<\mathbf{e}/q| \tau^-_0>\mathbf{e}/q)\,,\;\;\;t\ge0,\,\;\;\;x,y>0
\label{CSP}
\end{equation}
where $\mathbf{e}$ is an independent and exponentially distributed random variable with mean $1$.
It turns out that the measure on the left hand side can also be constructed as the result of
a Doob $h$-transform of $X$ killed when it
first exists $(0,\infty)$, i.e. at time
$\tau^{-}_0=\inf\{t>0:X_t\le0\}$. The resulting semi-group is thus given by
\begin{equation}
\p_x^{\uparrow}(X_{t}\in dy)=
\frac{h(y)}{h(x)}\p_x(X_t\in dy\,,t<\tau^{-}_0)
\,,\;\;\;t\ge0,\,\;\;\;x,y>0.
\label{doob}
\end{equation}
Here the function $h$ satisfies
 \[
 h(x)=\e\left(\int_0^\infty\ind_{\{I_t\ge -x\}}\,dL_t\right)
\,,\;\;\;x\ge0\, ,
\]
$L$ is the local time at zero of the reflected process $(X_t-\inf_{s\leq t} X_s, t\geq 0)$.
The family of measures $(\p^{\uparrow}_x, \, x>0)$ induced on $\mathcal{D}$ are probability measures 
if and only if $\psi'(0+)\leq 0$. 
When $X$ has unbounded variation paths,  the law $\p^\uparrow_x$ converges weakly as $x\downarrow 0$ to a measure denoted by $\p^\uparrow$. In the special case where $X$ oscillates, the function
$h$ satisfies $h(x)=x$. 

Now, define $\widehat{X}:=-X$, the dual process of $X$. Note that, under $\px$, the process $\widehat{X}$ is a L\'evy process with no positive jumps starting from the state $-x$. We denote by $\widehat{\p}_x$ for its law.
The dual process conditioned  to stay positive in the sense of (\ref{CSP}) is again a Doob $h$-transform of $(X, \widehat{\p}_x)$ killed when it
first exists $(0,\infty)$. In this case however one works with the $h$-function
$\widehat{h}(x)= \widehat{\e}\left(\int_0^\infty\ind_{\{I_t\ge -x\}}\,dL_t\right)$. This time, the resulting measure is a probability measure if and only if $\psi'(0+)\geq 0$. In the sequel, we shall work with a slightly different definition of $\widehat{\p}_x^\uparrow$ to the obvious analogue of $\p^\uparrow_x$ mentioned above.

Following the definition in Chapter VII of Bertoin \cite{Be}, for each $x>0$, the measure  $\widehat{\p}^\uparrow_x$ is defined as the result of  a Doob $h$-transform with respect to $\widehat{\p}_x$ of the kind (\ref{doob}) but with $h$-function given by $\widehat{h}(x)=W(x)$, where $W$ is the so-called scale function for the spectrally negative L\'evy process $\widehat{X}$. The latter is the
unique continuous increasing function with Laplace transform
\begin{equation}
\int_0^{\infty}e^{-\lambda x}W(x)\ud x=\frac{1}{\psi(\lambda)}, \qquad \lambda \geq 0.
\label{SF}
\end{equation}
In that case the measure $\widehat{\p}_x^\uparrow$ is always a probability measure. Note that $W(x)= \widehat{\e}\left(\int_0^\infty\ind_{\{I_t\ge -x\}}\,dL_t\right)$ precisely when $\psi'(0+)\geq 0$.  With this definition of  $\widehat{\p}_x^\uparrow$, there is always weak convergence as $x\downarrow 0$ to a probability measure which we denote by $\widehat{\p}^\uparrow$.

\subsection{Continuous state branching processes}

Continuous state branching processes are the analogue of Galton-Watson processes in continuous time and continuous state space. Such class of processes have  been introduced by Jirina \cite{ji} and studied by many authors included Bingham \cite{bi}, Grey\cite{gre}, Grimvall \cite{gri}, Lamperti \cite{la, la2}, to name but a few.
A continuous state branching process $Y=(Y_t, t\geq 0)$ is a Markov process taking values in $[0,\infty]$, where $0$ and $\infty$ are two absorbing states. Moreover, $Y$  satisfies the branching property; that is to say, the Laplace tranform of $Y_t$ satisfies
\begin{equation}
\ex( e^{-\lambda  Y_t})=\exp\{-x
u_t(\lambda)\},\qquad\textrm{for }\lambda\geq 0,
\label{CBut}
\end{equation}
for some function $u_t$. According to Silverstein \cite{si}, the function $u_t(\lambda)$ is determined by the
integral equation
\begin{equation}
\int_{u_t(\lambda)}^\lambda \frac{1}{\psi(u)}{\rm d} u=t
\label{DEut}
\end{equation}
where $\psi$ is the Laplace exponent of a spectrally positive L\'evy process.

Lamperti \cite{la1} observed that  continuous state branching
processes are connected to L\'evy processes with no negative jumps
by a simple time-change. More precisely, consider the spectrally positive L\'evy process
$(X,\px)$ started at $x>0$ and with Laplace exponent $\psi$. Now, we
introduce the clock
\[
A_{t}=\int_0^{t}\frac{\ud s}{X_s}, \qquad t\in [0,\tau_0^{-}).
\]
and its right-inverse $\theta(t)=\inf\{s\geq 0:\,A_s>t\}.$
Then, the time change process $Y=(X_{\theta(t)}, t\geq 0)$, under
$\px$, is a continuous state branching process
(or CB-process for short) with initial population of size $x$.
The transformation described above will henceforth be referred to as the {\it CB-Lamperti representation}.

In respective order, a CB-process is  called supercritical, critical or subcritical accordingly as its associated
L\'evy process drifts to $+\infty$, oscillates
or drifts to $-\infty$, in other words accordingly as $\psi'(0+)<0$, $\psi'(0+)=0$ or $\psi'(0+)>0$. It is known that if the CB-process $Y$ with branching mechanism $\psi$
satisfies that
\begin{equation}\label{cond}
\int_1^{\infty}\frac{\ud u}{\psi(u)}<\infty\, ,\end{equation} hence
$Y$ has a finite time extinction almost surely.

 In this work, we are also interested in CB-processes with immigration.
 In the remainder of this subsection, we assume that the
CB-process is critical, i.e. $\psi'(0+)=0$. Recall that a CB-process with immigration (or CBI-process) is a
strong Markov process taking values in $[0,\infty]$, where $0$ is no
longer absorbing. If $(Y^\uparrow_t : t\geq 0)$ is a process in this class, then its semi-group is  characterized by 
\[
 \ex(e^{-\lambda Y^\uparrow_t}) = \exp\{-x u_t(\lambda) - \int_0^t \phi(u_{t-s}(\lambda))ds\}\qquad\textrm{for }\lambda\geq 0,
\]
where $\phi$ is a Bernstein function satisifying $\phi(0)=0$ and is referred to as the immigration mechanism. See  for example Lambert \cite{lamb} for a formal
definition. Roelly and Rouault \cite{rr}, and more recently Lambert
\cite{lamb}, show that, if 
\[
T^-_0=\inf\{t>0: Y_t  =0 \},
\] 
then the limit
\begin{equation}
 \lim_{s\uparrow\infty}\px(Y_t \in dy| T^-_0>t+s)\,,\;\;\;t\ge0,\,\;\;\;x,y>0
\label{CBCSP}
\end{equation}
exists and defines a semi-group which is that of a CBI-process having initial population size $x$ and immigration mechanism
\[
\phi(\lambda)=\psi'(\lambda), \qquad \lambda\geq 0.
\]
The limit (\ref{CBCSP}) may be thought of as conditioning the CB-process to not to become extinct.

Lambert \cite{lamb} also proved an interesting connection between the conditioning (\ref{CBCSP}) for CB-process and (\ref{CSP}) for the underlying L\'evy process. Specifically he showed that
$(Y, \px^\uparrow) = (Y^\uparrow, \px)$ where the latter process has immigration mechanism given by $\psi'(\lambda)$. Another way of phrasing this is that the CBI-process obtained by conditioning a critical CB-process not to become extinct is equal in law to the underlying spectrally positive L\'evy process conditioned to stay positive an time changed with the CB-Lamperti representation.

\subsection{Stable processes and pssMp-Lamperti representation.}
Stable L\'evy processes with no negative jumps are L\'evy processes with Laplace exponent of the type (\ref{lk}) which satisfy the scaling property for some
index $\alpha>0$. More precisely, there exists a constant $\alpha>0$
such that for any $k>0$,
\begin{equation}
\textrm{ the law of  }\quad\left(kX_{k^{-\alpha}t}, t\geq 0\right)\quad \textrm{under}\quad   \mathbb{P}_x\,\textrm{ is }\, \mathbb{P}_{kx}.
\label{selfsim}
\end{equation}
In this subsection,  $(X, \mathbb{P}_x)$ will denote a stable L\'evy process with no negative jumps
of index $\alpha\in(1,2]$ starting at $x\in \R$,
(see Chapter VII in Bertoin \cite{Be} for further discussion on stable L\'evy processes). It is known, that the Laplace exponent of $(X, \mathbb{P}_x)$
takes the form
\begin{equation}
\psi(\lambda)=c_+\lambda^{\alpha} ,\qquad \lambda \geq 0,
\label{c+}
\end{equation}
where $c_+$ is a nonnegative constant.  Moreover, the density of 
its L\'evy measure is given by
$\Pi(\ud x)=c_+
x^{-\alpha-1}\ind_{\{x>0\}}\ud x$. The  case $\alpha= 2$ corresponds 
the process $(X, \mathbb{P}_x)$ being a multiple of 
 Brownian motion. In the remainder of this work, when we consider the case $\alpha=2$ we will refer to the Brownian motion, i.e. that we choose $c_+=1/2$.

Recall that the stable L\'evy process killed at the first time that it enters the negative half-line is defined by
\[
X_t\ind_{\{t<\tau^{-}_{0}\}}, \qquad t\geq 0,
\]
where $\tau^{-}_0=\inf\{t\geq 0:X_t\leq 0\}$. From the previous subsection, a stable L\'evy
process with no negative jumps conditioned to stay positive is tantamount to a Doob-$h$ transform of the killed process where $h(x)=x$. According to
Caballero and Chaumont \cite{cc}, both the process $X$ and its conditioned version belong to the class
of positive self-similar Markov processes; that is to say 
positive Markov processes satisfying  the property (\ref{selfsim}).

From Lamperti's work \cite{la3} it is known that the family of positive self-similar Markov processes up to
its first hitting time at $0$ may be expressed as the exponential of a L\'evy process, time changed by the inverse of its exponential functional.
More precisely, let $(X,\mathbb{Q}_x)$ be a self-similar Markov
process started from $x>0$  that fulfills the scaling property for
some $\alpha>0$, then under $\mathbb{Q}_x$ , there exists a L\'evy
process  $\xi=(\xi_{t}, t\geq 0)$ possibly killed at an independent exponential time which does not depend on $x$ and such that
\begin{equation}\label{lamp}
X_{t}=x\exp\Big\{\xi_{\zeta(tx^{-\alpha})}\Big\},\qquad 0\leq
t\leq x^{\alpha}I(\xi),
\end{equation}
where
\begin{equation*}
\zeta(t)=\inf\Big\{s\geq 0: I_{s}(\alpha\xi)>t\Big\},\quad
I_{s}(\alpha \xi)=\int_{0}^{s}\exp\big\{\alpha\xi_{u}\big\}\ud
u\textrm{ and } I(\alpha\xi)=\lim_{t\to+\infty}I_{t}(\alpha\xi).
\end{equation*}
We will refer to this transformation as {\it pssMp-Lamperti representation}.

In \cite{ckp},  it was shown that the
Laplace exponent of the underlying L\'evy process, $\xi$, of
 the stable L\'evy process with index $\alpha\in(1,2)$ killed at the first time that it 
enters the negative half-line which  is
given by
\begin{equation}\label{chalev}
\Psi(\lambda)=m\frac{\Gamma(\lambda+\alpha)}{\Gamma(\lambda)\Gamma(\alpha)},\qquad \textrm{ for }\quad \lambda\geq 0,
\end{equation}
where $m>0$ is the mean of $-\xi$ which is finite. Note that this last fact implies that the process $\xi$ drifts towards $-\infty$.  In the Brownian case, i.e when $\alpha=2$, we have that the L\'evy process $\xi$ is  a Brownian motion with drift $a=-1/2$.
 
The Laplace
exponent  of the underlying L\'evy process, denoted by
$\xi^*$, of the 
stable L\'evy process  (with $\alpha\in(1,2)$) conditioned to stay positive is also computed in \cite{ckp}. It is given by 
\begin{equation}
\Psi^*(\lambda)=m^*\frac{\Gamma(\lambda-1+\alpha)}{\Gamma(\lambda-1)\Gamma(\alpha)}\qquad \textrm{ for }\quad \lambda\geq 0,
\end{equation}
where $m^*>0$ is the mean of $\xi^*$ which is also finite. In this case, the L\'evy process $\xi^*$ drifts towards $+\infty$.  When $\alpha=2$,  it is not difficult to show that the process $\xi^*$ is a Brownian motion with drift $a=1/2$.\\
Finally, we know that under $\widehat{\p}_x$,
the stable L\'evy process $X$ has no positive jumps. In this case, the density of 
its L\'evy measure is given by
$\Pi(\ud x)=c_{-}
|x|^{-\alpha-1}\ind_{\{x<0\}}\ud x,$ where $c_{-}=c_+$.
From  Corollary 6 in \cite{ckp}, it is known that the underlying Levy process in the pssMp-Lamperti representation of 
the spectrally negative stable L\'evy process conditioned to stay positive  is $\widehat{\xi}$, the dual of $\xi$. Note that in the case $\alpha=2$, the processes $\widehat{\xi}$ and $\xi^*$  are the same. 

It is not difficult to show that the process $\xi$ corresponds to $\xi^*$ conditioned to drift towards $-\infty$ (or equivalently $\xi^*$ is $\xi$ conditioned to drift to $+\infty$). This relation will be used later and hence  we register it as a proposition below, its proof can be found in \cite{ckp}. In the sequel, $P$ will be a reference probability measure on $\mathcal{D}$ under which $\xi$ and $\xi^*$  are L\'evy processes whose repective laws are defined above. 
\begin{proposition}\label{prop1} For every, $t\geq 0$, and every bounded measurable function $f$, 
\[
E\big(f(\xi^*_t)\big)=E\big(\exp\{\xi_t\}f(\xi_t)\big).
\]
In particular, the process $-\xi^*$ and $\xi$ satisfy the Cram\'er's condition, i.e.
\[
E\big(\exp\{-\xi^*_1\}\big)=1\qquad\textrm{and}\qquad E\big(\exp\{\xi_1\}\big)=1.
\]
\end{proposition} 
\section{Time reversed CB-processes and total progeny.}
In this section, we shall dwell on two aspects of the paths of general CB-processes which shall be developed in more detail for the specific case that the underlying L\`evy process is $\alpha$-stable in later sections.
Specifically, we
are interested in  time reversal properties and distributional features of the total progeny of    CB-processes. 

First, we will determine the law of the time reversed process
$(Y_{(T_0^- -t) -}, \,0\leq t< T_0^-)$ under $\{T_0^-<\infty\}$, where $T_0^-$ is the extinction time, i.e. $T_0^-=\inf\{t:Y_t=0\}$. Our arguments are based on the CB-Lamperti representation, path decomposition and  time reversal properties of L\'evy processes with no negative jumps.

One motivation for this time reversed analysis is the recent work of Krikun \cite{kri}, where it is showed that the local growth of random quadrangulations is governed by certain critical time reversed branching process and moreover that its rescaled profile converge in the sense of Skorokhod to the time reversed stable CB-process with index $\alpha=3/2$, see  Theorem 4 in \cite{kri}.\\

 The total progeny until time $t\geq 0$, of a CB-processes  is defined as follows
\[
J_t=\int_{0}^{t}Y_u\ud u.
\]
At the end of this section, we will provide distributional identities for $J_{T^+_a}$ where
\[
T^+_a=\inf\{t>0: Y_t\geq a\}.
\]

\subsection{Time reversal.}

Let $\tau^-_x = \inf\{t>0 : X_t \leq x\}$ and $\tau^+_x = \inf\{t>0 : X_t \geq x\}$, be the first passage time of $X$ below and above $x\in \R$, respectively. We also introduce
\[
 \sigma_x = \sup\{t>0 : X_t \leq x\},
 \]
the last passage time of $X$ below  $x\in \R$. For the associated CB-process $Y$, we define
$ T_x^- = \inf\{t>0 : Y_t \leq x\}$,
its first passage time  below  $x\in \R_+$. 

Now for $x>0$, let  $Z^{(x)}=(X_{\sigma_x +t}, \, t\geq 0)$ and denote by $\theta^{(x)}$ for the right-inverse of the functional
\[
A^{(x)}=\int_{0}^{t}\frac{\ud s}{Z^{(x)}_{s}}\quad\textrm{for } \quad\, t\geq 0
\]
Recall the notation  $\widehat{\mathbb{P}}^\uparrow$ for the law of $\widehat{X}$ conditioned to stay positive in the sense of the previous section.
We remark that, under  $\widehat{\mathbb{P}}^\uparrow$, the canonical process $X$ drifts towards $\infty$ and also that $X_t>0$ for $t>0$.

\begin{proposition}\label{prop2} 
Let $y>0$. Then for every $x\in(0,y)$, the law of the time reversed process $(Y_{(T_x^- -t) -}, \,0\leq t< T_x^-)$ under
$\mathbb{P}_y(\cdot|\tau^{-}_x<\infty)$ is the same
as that of the shifted process
$(Z^{(x)}_{\theta^{(x)}(t)},0\leq t< A^{(x)}_{\sigma_{y}})$ under $\widehat{\p}^{\uparrow}$.
\end{proposition}
\noindent {\it Proof:}
By Theorem VII.18 of Bertoin \cite{Be}, we know that for $y>0$
\begin{equation}\label{B0}
\Big\{(X_t,\,0\leq t< \sigma_{y} ), \widehat{\mathbb{P}}^\uparrow\Big\}\stackrel{d}{=}
\Big\{(X_{(\tau^-_0-t)-}, \, 0\leq t<\tau^-_0), \mathbb{P}_{y}(\cdot|\tau^-_{0}<\infty)\Big\},
\end{equation}
which implies that for $x\in(0,y)$,
\begin{equation}
\Big\{(X_t,\,\sigma_x\leq t< \sigma_{y} ), \widehat{\mathbb{P}}^\uparrow\Big\}\stackrel{d}{=}
\Big\{(X_{(\tau^-_x-t)-}, \, 0\leq t<\tau^-_x), \mathbb{P}_{y}(\cdot|\tau^-_{x}<\infty)\Big\}.
\label{B}
\end{equation}
Next from the definition of $Y$, under $\p_y(\cdot|\tau^-_{x}<\infty)$, we have
\begin{equation}
(Y_{(T_x^--t)-},\, 0\leq t <T_x^-)=(X_{\theta(A_{\tau^-_x} -t) -}, 0\leq t< A_{\tau^-_x}).
\label{CSB}
\end{equation}
Define
\[
\theta'(t) = \inf\{s>0 : B_s >t\}\quad
\text{ where }\quad
B_s = \int_0^s \frac{1}{X_{\tau^-_{x} -u} }du.
\]
Setting $t= B_s$,  we have
\[
A_{\tau^-_x} - B_s = \int_0^{\tau^-_x}\frac{1}{X_u} du -  \int_0^s \frac{1}{X_{\tau^-_{x} -u} }du =
\int_0^{\tau^-_x -s}\frac{1}{X_u}du
\]
and hence
\[
X_{\theta(A_{\tau^-_x} - t)-} = X_{\theta(A_{\tau^-_x - s})-} = X_{(\tau^-_x - s)-} = X_{(\tau^-_x - \theta'(t))-}.
\]
Note also that $T_x^- = A_{\tau^-_x} = B_{\tau^-_x}$. From (\ref{B}) the latter, under $\mathbb{P}_y(\cdot|\tau^-_{x}<\infty)$, is equal
in distribution to
\[
\int_{\sigma_x}^{\sigma_y}\frac{\ud t}{X_t}\qquad \textrm{under }\quad \widehat{\mathbb{P}}^\uparrow.
\]
Now,
it follows from (\ref{B})  and (\ref{CSB}) that
\[
\Big\{(Y_{(T_x^--t)-},\,0\leq t< T_x^-), \mathbb{P}_y(\cdot|\tau^-_x <\infty)\Big\} \stackrel{d}{=}\Big\{\big(Z^{(x)}_{\theta^{(x)}(t) },\, 0\leq t< A^{(x)}_{\sigma_y}\big), \widehat{\mathbb{P}}^\uparrow\Big\}
\]
as required.\QED

\begin{theorem}\label{trcb-thrm} 
If  condition (\ref{cond}) is satisfied, then for every $y>0$
\[
\Big\{(Y_{(T_0^- -t) -}, \, 0\leq t< T^-_0 ), \mathbb{P}_y\Big\} \stackrel{d}{=}
\Big\{(X_{\theta(t) },\,0\leq t< A_{\sigma_y}  ), \widehat{\mathbb{P}}^\uparrow\Big\},
\]
where $\stackrel{d}{=}$ denotes equality in law or distribution.
\end{theorem}
\noindent{\it Proof:} We first prove that the CB-Lamperti representation is well defined for the process $(X,\widehat{\p}^{\uparrow})$. In order to do so, 
it is enough to prove that the map $s\mapsto 1/X_s$ is integrable in a neighbourhood of 0.  Take $\epsilon>0$  small enough  and note that 
\[
\int_{0+}\frac{\ud s}{X_s}\ind_{\{X_s\leq \epsilon\}}<\infty \qquad \textrm{if and only if}\qquad \int_{0+}\frac{\ud s}{X_s}<\infty
\]
Now, from Proposition VII.15 in \cite{Be} and Fubini's Theorem, we get that for $t>0$
\begin{equation}\label{riv1}
\widehat{\e}^{\uparrow}\left(\int_0^t\frac{\ud s}{X_s}\ind_{\{X_s\leq \epsilon\}}\right)=k\int_0^{t} n\left(\frac{W(X_s)}{X_s}\ind_{\{X_s\leq \epsilon\}};\,s<\zeta\right)\ud s,
\end{equation}
where $W$ is the scale function defined in (\ref{SF}), $k$ is a strictly positive constant which only depends on the normalization of local time 
$L$ at zero of the reflected process $(X_t-\inf_{s\leq t}X_s,\,t\geq 0)$, $n$ its excursion measure and $\zeta$ denotes the life time of 
the generic excursion.

The identity (\ref{riv1}) implies that
\[
\begin{split}
\widehat{\e}^{\uparrow}\left(\int_0^t\frac{\ud s}{X_s}\ind_{\{X_s\leq \epsilon\}}\right)&\leq k\int_0^{\infty} n\left(\frac{W(X_s)}{X_s}\ind_{\{X_s\leq \epsilon\}};\,s<\zeta\right)\ud s,\\
&=k \,n\left(\int_0^{\zeta}\frac{W(X_s)}{X_s}\ind_{\{X_s\leq \epsilon\}}\ud s\right).
\end{split}
\]

On the one hand, from the occupation measure of the excursion law, we deduce that
 \[
 \,n\left(\int_0^{\zeta}\frac{W(X_s)}{X_s}\ind_{\{X_s\leq \epsilon\}}\ud s\right)=\int_{0}^{\infty} \frac{W(x)}{x}\ind_{\{x\leq \epsilon\}}\ud \mathcal{V}(x),
 \]
 where $\ud \mathcal{V}$ is the renewal measure of $H=(H_t,t\geq 0)$, the upward ladder height process of $(X,\widehat{\p})$ 
(see \cite{Be} for a proper definition). On the other hand, in this particular case the process $H$ is a pure drift which 
implies that its renewal measure $\ud \mathcal{V}$ is in fact, the Lebesgue measure. Hence, we get
\[
\widehat{\e}^{\uparrow}\left(\int_0^t\frac{\ud s}{X_s}\ind_{\{X_s\leq \epsilon\}}\right)\leq k\int_{0}^{\epsilon} \frac{W(x)}{x}\ud x.
\] 
  
It is known (see for instance the proof of Proposition VII.10 in \cite{Be}) that there exist two positive constants $0<c_1<c_2$ such that
\[
c_1\frac{1}{x\psi(1/x)}\leq W(x)\leq c_2\frac{1}{x\psi(1/x)} \qquad \textrm{ as }\quad x\to 0.
\]
Hence, we have
\[
\int_{0}^{\epsilon} \frac{W(x)}{x}\ud x\leq c_2\int_{1/\epsilon}^{\infty}\frac{\ud u}{\psi(u)},
\]
which is finite from our hypothesis. We may now conclude that the map $s\mapsto 1/X_s$ is integrable on a neighbourhood of 0.\\
Now, we may follow the proof Proposition \ref{prop2} line by line replacing $T_x^-$ by $T_0^-$ and  get the desired result.\QED
\begin{corollary}
Suppose that  the L\'evy process $(X, \p)$
does not drift towards $+\infty$. Then  for every $x>0$ and $0<y\leq x$,
\[
\p_x\Big(\inf_{0\leq t\leq U_y}Y_t\geq z\Big)=\frac{W(y-z)}{W(y)}\ind_{\{z\leq y\}},
\]
where $U_y=\sup\{t>0:Y_t\geq y\}$ and the scale function $W$ satisfies (\ref{SF}).
\end{corollary}
\bigskip \noindent{\it Proof:} From Proposition \ref{prop2} and since $X$ has no negative jumps, it is clear that
\begin{equation}\label{trcb}
\Big\{(Y_{(T_0^- -t) -}, \, T_0^- -U_y\leq t\leq T  ), \mathbb{P}_x)\Big\} \stackrel{d}{=}
\Big\{(X_{\theta(t) },\,A_{\tau^+_{y}}\leq t< A_{\sigma_x}  ), \widehat{\mathbb{P}}^\uparrow\Big\}.
\end{equation}
On the other hand, by Theorem 1 in \cite{CD}, we have that for $z\leq y$
\[
\widehat{\mathbb{P}}_y^{\uparrow}\Big(\inf_{t\geq 0 } X_t\geq z\Big)=\widehat{\mathbb{P}}_y^{\uparrow}\Big(\inf_{0\leq t\leq \sigma_x } X_t\geq z\Big)=\frac{W(y-z)}{W(y)}.
\]
Hence from (\ref{trcb}), the above formula and the Markov property of $(X,\widehat{\p}^{\uparrow})$, the  statement of the corollary follows.\QED
\subsection{Total progeny.}
It is known that  $\tau^-=(\tau^-_{-x}, x\geq 0)$, the first passage time process of $X$, is a (possibly killed)    subordinator with Laplace exponent $\Phi(q)$.  The killing rate is given by $\Phi(0)$ and recall that
$X$ drifts to $+\infty$ if and only if $\Phi(0)>0$.

Now for $q\geq 0$, we define the scale functions $W^{(q)}$ and $Z^{(q)}$, both mapping $\mathbb{R}$ to $[0,\infty)$,  as follows. In the first case, $W^{(q)}(x) = 0 $ for $x<0$ and otherwise it is the unique continuous function with Laplace transform
\[
\int_0^{\infty}e^{-\lambda x}W^{(q)}(x)\ud x=\frac{1}{\psi(\lambda)-q}\qquad \textrm{for} \quad \lambda>\Phi(q).
\]
In the second case, for $x\in\mathbb{R}$ 
\[
Z^{(q)}(x)=1+q\int_{0}^x W^{(q)}(y)\ud y.
\] Note that from the definition  of 
$W^{(q)}$, we have that $W^{(0)}=W$, which was defined in Section 2.

\begin{theorem} Let $(X, \px)$ be a L\'evy process with no negative jumps starting from $x$ and $(Y, \px)$ its associated CB-process. Then,
\begin{itemize}
\item[i)] For each $a\geq x>0$ and $q\geq 0$,
\[
\e_x\left(\exp\left\{-q\int_{0}^{T^+_a}Y_s\ud s\right\}\ind_{\{T^+_a<T\}}\right)=Z^{(q)}(a-x)-W^{(q)}(a-x)\frac{Z^{(q)}(a)}{W^{(q)}(a)}.
\]
\item[ii)] For each $a\geq x>0$ and $q\geq 0$,
\[
\e_x\left(\exp\left\{-q\int_{0}^{T_0^-}Y_s\ud s\right\}\ind_{\{T_0^-<T^+_a\}}\right)=\frac{W^{(q)}(a-x)}{W^{(q)}(a)}.
\]
\end{itemize}
\end{theorem}
\bigskip \noindent{\it Proof:} From the CB-Lamperti representation under $\px$, we have that
\[
\tau^{+}_a=\int_{0}^{T^+_a}Y_s\ud s\quad\textrm{ and }\quad \tau^-_0=\int_{0}^{T_0^-}Y_s\ud s.
\]
Now, the result follows  from an application of Theorem $8.1$ in \cite{ky} for the L\'evy process $X$. Recall that the process $X$ has no negative jumps and to implement the aforementioned result, which applies to spectrally negative processes, one must consider the problem of two-sided exit from $[0,a]$ of $-X$ when $X_0=a-x$.\QED

\section{Self-similar CB-processes}
Suppose now  that $(X,\px)$ is a spectrally positive $\alpha$-stable process with index $\alpha \in(1,2]$ starting from $x>0$. We refer to $Y$, the associated continuous state branching process, as the $\alpha$-stable CB-process.

In this section, we are interested in the self-similar property of the $\alpha$-stable CB-process. Since, it is is a positive  Markov process we will determine its underlying L\'evy process in the pssMp-Lamperti representation. Such representation will be important to study its asymptotic behaviour at the extinction time $T_0^-$.    

We start with a generic result which, in some sense, is well known folk-law  and will be useful throughout the remainder of this section. For the sake of completeness we include its proof.

\begin{proposition}\label{generic}
Suppose that $X$ is any positive self similar Markov process with self-similarity index $\alpha$ and let $\theta$ be the CB-Lamperti time change. Then $X_\theta$ is a positive self similar Markov process with self-similarity index $\alpha-1$ with the same underlying L\'evy process as $X$.
\end{proposition}

\bigskip \noindent{\it Proof:} Suppose that $\eta$ is the underlying L\'evy process for the process $X$. We first define,
\[
A_t=\int_0^t\frac{\ud s}{X_s},\qquad I_t(\alpha\eta)=\int_0^te^{\alpha\eta_s}\ud s \qquad \textrm{and}\qquad I_t\big((\alpha-1)\eta\big)=\int_0^te^{(\alpha-1)\eta_s}\ud s.
\]
Recall that $\zeta$ is the right-continuous inverse of $I(\alpha\eta)$. From the pssMp-Lamperti transform of $X$ and the change of variable $s=x^{\alpha}I_u(\alpha\eta)$, we get that 
\[
A_{x^{\alpha}I_t(\alpha\eta)}=\int_{0}^{x^{\alpha}I_t(\alpha\eta)}\frac{\ud s}{x\exp\{\eta_{\zeta(s/x^{\alpha})}\}}=x^{\alpha-1}\int_0^{t}\frac{e^{\alpha\eta_u}}{e^{\eta_u}}\ud u=x^{\alpha-1}\int_0^te^{(\alpha-1)\eta_u}\ud u.
\] 
On the other  hand, the right-continuous inverse of $I\big((\alpha-1)\eta\big)$ is defined by 
\[
h(t)=\inf\Big\{s\geq 0: I_s\big((\alpha-1)\eta\big)>t\Big\},
\]
and recall that  $\theta$ is the right-continuous inverse function of $A$. Hence,  we have that for any $0\leq t< x^{\alpha-1}I_{\infty}\big((\alpha-1)\eta\big)$,
\[
\begin{split}
h(t/x^{\alpha-1})&=\inf\Big\{s\geq 0: I_s\big((\alpha-1)\eta\big)>t/x^{\alpha-1}\Big\}=\inf\Big\{s\geq 0: A_{x^{\alpha}I_u(\alpha\eta)}>t\Big\}\\
&=\inf\Big\{\zeta(u/x^{\alpha})\geq 0:A_{u}>t\Big\}=\zeta(\theta(t)/x^{\alpha}).
\end{split}
\]
From the pssMp-Lamperti representation of $X$, we have that, under $\px$, for every $0\leq t<T_0^-$,
\[
X_{\theta(t)}=x\exp\Big\{\eta_{\zeta(\theta(t)/x^{\alpha})}\Big\}=x\exp\Big\{\eta_{h(t/x^{\alpha-1})}\Big\}
\]
thus completing the proof.\QED

\subsection{The pssMp-Lamperti representation of the $\alpha$-stable CB-process.}
Recall  that  there exists a spectrally positive L\'evy process, $\xi=(\xi_t,t \geq 0)$ starting from $0$, which drifts towards $-\infty$,  whose Laplace exponent is given by
\[
 \Psi(\theta) = m\frac{\Gamma(\theta+\alpha)}{\Gamma(\theta)\Gamma(\alpha)},
\]
and does not depend on $x$. When $\alpha=2$, we recall that $\xi_t=B_t-t/2$, where $B$ is a standard Brownian motion.
We have the following direct corollary to Proposition \ref{generic}.

\begin{corollary}\label{noarrow}
The process $Y$ is a positive self similar  Markov process with self-similarity index $\alpha-1$. Moreover,  its pssMp-Lamperti representation  under $\mathbb{P}_x$ is given by
\[
Y_t=x\exp\Big\{ \xi_{h(tx^{-(\alpha-1)})}\Big\}, \qquad 0\leq t\leq x^{\alpha-1}\int_0^{\infty}\exp\Big\{(\alpha-1)\xi_u\Big\}\ud u,
\]
 where
\[
h(t) = \inf\left\{ s\geq 0 : \int_0^s \exp{\Big\{(\alpha -1)\xi_u \Big\}}du > t \right\}.
\]
\end{corollary}

From the above representation of $Y$, we get that $T_0^-=x^{\alpha-1}I_{\infty}\big((\alpha-1)\xi\big)$.
As a prelude to the next theorem we shall first prove the following auxiliary lemma which says in particular that the distribution of $I_{\infty}\big((\alpha-1)\xi\big)$ has a Fr\'echet distribution and moreover that the Fr\'echet distribution is self-decomposable.
\begin{lemma}
The distribution of $I:=I_{\infty}\big((\alpha-1)\xi\big)$ is given by 
\begin{equation}\label{felp}
P(I \leq t) = \exp\{-[c_+(\alpha - 1) t]^{-1/(\alpha -1)}\}.
\end{equation}
Moreover, $I$ is self-decomposable and it has a completely monotone density with respect to the Lebesgue measure and is given by
\begin{equation}\label{density}
P(I \in \ud x) = c_+(\alpha-1)\big(c_+(\alpha-1)x\big)^{-\alpha/(\alpha-1)}\exp\{-[c_+(\alpha - 1) x]^{-1/(\alpha -1)}\}\ind_{(0,\infty)}\ud x.
\end{equation}
\end{lemma}
\bigskip \noindent{\it Proof:} From the pssMp-Lamperti representation of $(Y,\px)$, we deduce that $T^-_0=x^{\alpha}I$. From Bingham \cite{bi}, it is known that
\[
\mathbb{P}_x(T \leq t) = e^{- xu_t(\infty)}
\]
where $u_t(\infty)$ solves
\[
\int_{u_t(\infty)}^\infty \frac{1}{c_+v^\alpha}dv = t.
\]
Therefore
\[
\mathbb{P}_x(T^-_0 \leq t)=P(I \leq t/x^{\alpha-1}) = \exp\{-x[c_+(\alpha - 1) t]^{-1/(\alpha -1)}\},
\]
which implies (\ref{felp}). 

Let $a<0$, then
\[
I=\int_0^{\infty}e^{(\alpha-1)\xi_u}\ud u=\int_0^{S^{-}_a}e^{(\alpha-1)\xi_u}\ud u+e^{(\alpha-1)a}\int_0^{\infty}e^{(\alpha-1)\xi'_u}\ud u,
\]
where $\xi'=(\xi_{S^{-}_a+t}-a, t\geq 0)$ and $S^{-}_a=\inf\{t\geq 0: \xi_t\leq a\}$. Then, self-decomposability follows from the independence of  $(\xi_{t},0\leq t\leq \tau^{-}_a)$ and $\xi'$. Self-decomposable distributions on $\R_+$ are unimodal (see for instance Chapter 10 in Sato \cite{Sa}), i.e. that they have a completely monotone density on $(0,\infty)$, with respect to the Lebesgue measure.

\begin{theorem}\label{Trscb}
For each $x>0$
\begin{equation}\label{trscb}
\Big\{(Y_{(T_0^- -t) -}  : t< T^-_0 ), \mathbb{P}_x\Big\} \stackrel{d}{=}
\Big\{(X_{\theta(t)}, 0\leq t<A_{\sigma_x}),\widehat{\mathbb{P}}^\uparrow\Big\}.
\end{equation}
Moreover the process $(X_{\theta_{t}}, t\geq 0)$, under $\widehat{\p}^{\uparrow}$ is a positive self-similar Markov process with 
index $\alpha-1$, starting from $0$, with the same semigroup as the  processes $(X_\theta, \mathbb{P}^\uparrow_y)$  for $y>0$, 
and 
with entrance law given by
\begin{equation}\label{entlaw}
\widehat{\e}^{\uparrow}(f(X_{\theta(t)}))=\frac{c_+}{m}\int_0^{\infty}x^{-(2\alpha-1)/(\alpha-1)}f(tc_+(\alpha-1)/x)e^{-x^{-1/(\alpha-1)}}\ud x
\end{equation} 
where $t>0$ and $f$ is a positive measurable function.
In particular, under $\widehat{\p}^{\uparrow}_y$ for $y>0$ we have
\[
X_{\theta_s}=y\exp\Big\{ \widehat{\xi}_{\widehat{h}(s y^{-(\alpha-1)})}\Big\}, \qquad 0\leq s\leq y^{\alpha-1}\int_0^{\infty}\exp\Big\{(\alpha-1)\widehat{\xi}_u\Big\}\ud u,
\]
 where
\[
\widehat{h}(s) = \inf\left\{ u\geq 0 : \int_0^u \exp{\Big\{(\alpha -1)\widehat{\xi}_u \Big\}}du > s \right\}.
\] 
\end{theorem}
\bigskip \noindent{\it Proof:} The time reversal property follows from Theorem \ref{trcb-thrm}. The pssMp-Lamperti representation of the process $(X_{\theta(t)}, t\geq 0)$ under $\widehat{\p}^\uparrow_y$ when issued from $y>0$
 follows from Proposition \ref{generic}, noting in particular that ($X,\widehat{\mathbb{P} }^\uparrow_y) $ is a spectrally positive stable process conditioned to stay positive which a positive self-similar process with index $\alpha$.

The L\'evy process $\widehat{\xi}$ satisfies the conditions of Theorems 1 and 2 in \cite{cc1}. Hence, the family of processes $(X_{\theta(t)}, t\geq 0)$
  under $\widehat{\p}_y^\uparrow$, for $y>0$, converges weakly with respect to the Skorohod topology, as $y\downarrow 0$, towards a pssMp starting from $0$ which is 
$(X_{\theta(t)}, t\geq 0)$ under $\widehat{\mathbb{P}}^\uparrow$.  

From Proposition $3$ in \cite{cc1}, the entrance law of $(X_\theta,\widehat{\p}^\uparrow )$ is given by
\[
\widehat{\e}^{\uparrow}\Big(f(X_{\theta(t)})\Big)=\frac{1}{(\alpha-1)m}E\Big(I^{-1}f\big(tI^{-1}\big)\Big),
\]
for every $t>0$ and every $f$ positive and measurable function. Therefore, from the form of the density of $I$ given in the previous Lemma and some basic calculations, we get (\ref{entlaw}).\QED 
It is important to note that when $\alpha=2$, the process $(X_{\theta(t)},t\geq 0)$ under $\widehat{\p}^{\uparrow}$ is in fact the CB-process with immigration. 
This follows from the remark made in subsection 1.4 that in this particular case, we have that $\widehat{\xi}=\xi^*$. 

\subsection{Asymptotic behaviour at $T_0^-$}
We start by stating the integral test for the lower envelope of $(Y_{(T_0^- -t)^-}, 0\leq t\leq T_0^-)$, under $\px$, at $0$.
\begin{theorem}Let $f$ be an increasing fucntion such that $\lim_{t\to 0}f(t)/t=0$, then for every $x>0$
\[
\px\Big(Y_{(T_0^- -t)^-}<f(t), \textrm{ i.o., as } t\to 0\Big)= 0 \textrm{ or } 1, 
\]
accordingly as
\[
\int_{0+} f(t)\, t^{-\alpha/(\alpha-1)}\,\ud t\qquad \textrm{is finite or infinite.}
\]
In particular, 
\[
\liminf_{t\to 0}\frac{Y_{(T_0^- -t)^-}}{t^{\kappa}}=
\left\{ \begin{array}{ll}
0 &\textrm{if $\kappa>\frac{1}{\alpha-1}$}\\
+\infty &\textrm {if $\kappa\leq \frac{1}{\alpha-1}$}
\end{array} \right.
\qquad \p_x-\textrm{a.s.}
 \]
\end{theorem}
\bigskip \noindent{\it Proof:} Fix $x>0$.  From Theorem \ref{Trscb}, we deduce that
\[
\px\Big(Y_{(T_0^- -t)^-}<f(t), \textrm{ i.o., as } t\to 0\Big)=\widehat{\p}^{\uparrow}\Big(Z_t<f(t), \textrm{ i.o., as } t\to 0\Big).
\]
Since $(Y, \widehat{\p}^{\uparrow})$ is a pssMp with index $\alpha-1$ starting from $0$, from Theorem 3 in \cite{cp} it is enough to show that 
\begin{equation}\label{asy1}
P(I>t)\sim \big(c_+(\alpha-1)t\big)^{-1/(\alpha-1)}\qquad \textrm{as }t \textrm{ goes to }+\infty,
\end{equation}
to obtain the result.

From Lemma 1 and with the change of variable $h=\big(c_+(\alpha-1)t\big)^{-1/(\alpha-1)}$, we have that
\[
\lim_{t\to +\infty}\frac{P(I> t)}{\big(c_+(\alpha-1)t\big)^{-1/(\alpha-1)}}=\lim_{h\to 0}\frac{1-e^{-h}}{h}=1,
\] 
which proves (\ref{asy1}).\QED
Let us introduce $\mathcal{H}_0$ the class of increasing functions $f:(0,+\infty)\mapsto [0,+\infty)$
\begin{itemize}
\item[i)] $f(0)=0$ and
\item[ii)] there exist $\beta\in(0,1)$ such that $\displaystyle\sup_{t<\beta}\displaystyle\frac{t}{f(t)}<\infty$.
\end{itemize}
We also introduce the future infimum process of  
$(Z_t, \widehat{\p}^{\uparrow})$,
\[
J_t=\inf_{s\geq t} Z_s, \qquad \textrm{for all}\quad t\geq 0,  
\] 
and denote by $\underline{Y}_t$ for the infimum of the CB process $(Y,\p_x)$ over $[0,t]$.

The upper envelope of the process $(\underline{Y}_{(T_0^- -t)^{-}}, 0\leq t\leq T_0^-)$, under $\px$, at $0$ is described by the integral test in the following theorem.
\begin{theorem} 
Let $f\in \mathcal{H}_0$, then for every $x>0$
\begin{itemize}
\item[i)] If 
\[
\int_{0+}  \exp\Big\{-\big(c_+(\alpha-1)t/f(t)\big)^{-1/(\alpha-1)}\Big\}\frac{\ud t}{t}<\infty,
\]
then for all $\epsilon>0$
\[
\px\Big(\underline{Y}_{(T_0^- -t)^-}>(1+\epsilon)\big(f(t)\big)^{1/\alpha-1}, \textrm{ i.o., as } t\to 0\Big)= 0.
\]
\item[ii)] If 
\[
\int_{0+}  \exp\Big\{-\big(c_+(\alpha-1)t/f(t)\big)^{-1/(\alpha-1)}\Big\}\frac{\ud t}{t}=\infty,
\]
then for all $\epsilon>0$
\[
\px\Big(\underline{Y}_{(T_0^- -t)^-}>(1-\epsilon)\big(f(t)\big)^{1/\alpha-1}, \textrm{ i.o., as } t\to 0\Big)= 1.
\]
\end{itemize}
In particular, we have the following law of the iterated logarithm
\[
\limsup_{t\to 0}\frac{\underline{Y}_{(T_0^- -t)^{-}}}{t^{1/(\alpha-1)}(\log \log(1/t))^{1-\alpha}}=\big(c_+(\alpha-1)\big)^{1/(\alpha-1)}, \qquad \px-\textrm{a.s.}
\]
\end{theorem}
\bigskip \noindent{\it Proof:} Here, we will apply Theorem 1 in \cite{pa0}. First, we note again that from Theorem \ref{Trscb}, we have the following equality
\[
\begin{split}
\px\Big(Y_{(T_0^- -t)^-}>&(1+\epsilon)\big(f(t)\big)^{1/\alpha-1}, \textrm{ i.o., as } t\to 0\Big)\\
&=\widehat{\p}^{\uparrow}\Big(Z_t>(1+\epsilon)\big(f(t)\big)^{1/\alpha-1}, \textrm{ i.o., as } t\to 0\Big).
\end{split}
\] 

Hence, according to part $i)$ in Theorem 1 in \cite{pa0} and noting that the process $(Z,\widehat{\p}^{\uparrow})$ has no positive jumps, the right-hand side of the above equality is equal $0$, for all $\epsilon>0$, if
\[
\int_{0+} \exp\Big\{-\big(c_{+}(\alpha-1)t/f(t)\big)^{-1/(\alpha-1)}\Big\}\frac{\ud t}{t}<\infty.
\]
In order to prove part $ii)$, we note that from Theorem \ref{Trscb}, we have
\[
\begin{split}
\px\Big(Y_{(T_0^- -t)^-}>&(1-\epsilon)\big(f(t)\big)^{1/\alpha-1}, \textrm{ i.o., as } t\to 0\Big)\\
&=\widehat{\p}^{\uparrow}\Big(Z_t>(1-\epsilon)\big(f(t)\big)^{1/\alpha-1}, \textrm{ i.o., as } t\to 0\Big).
\end{split}
\]
Hence applying part $ii)$ of Theorem 1 in \cite{pa0}, we obtain that the above probability  is equal $1$, for all $\epsilon>0$, if
\[
\int_{0+}  \exp\Big\{-\big(c_{+}(\alpha-1)t/f(t)\big)^{-1/(\alpha-1)}\Big\}\frac{\ud t}{t}=\infty,
\]
and now the proof is complete.\QED
The following result describe the upper envelope of the time reversed processes  $(Y_{(T_0^- -t)^{-}}, 0\leq t\leq T_0^-)$ and $((Y-\underline{Y})_{(T_0^- -t)^{-}}, 0\leq t\leq T_0^-),$ by  laws of the iterated logarithm.
\begin{theorem} For every $x>0$, we have
\begin{equation}\label{lil1}
\limsup_{t\to 0}\frac{Y_{(T_0^- -t)^{-}}}{t^{1/(\alpha-1)}(\log \log(1/t))^{1-\alpha}}=\big(c_+(\alpha-1)\big)^{1/(\alpha-1)}, \qquad \px-\textrm{a.s.,}
\end{equation}
and
\begin{equation}\label{lil2}
\limsup_{t\to 0}\frac{(Y-\underline{Y})_{(T_0^- -t)^{-}}}{t^{1/(\alpha-1)}(\log \log(1/t))^{1-\alpha}}=\big(c_+(\alpha-1)\big)^{1/(\alpha-1)}, \qquad \px-\textrm{a.s.}
\end{equation}
\end{theorem}
\noindent{\it Proof:}  From Lemma 1, we have
\[
-\log P(I\leq t)=\Big(c_+(\alpha-1)t\Big)^{-1/(\alpha-1)}.
\]
Hence in this case, the condition of Theorem 6 in \cite{pa} is satisfied and (\ref{lil1}) follows.  Now, since the time reversed process $(Y_{(T_0^- -t)^{-}}, 0\leq t\leq T_0^-)$, under $\px$, is a positive self-similar Markov process starting from $0$, with no positive jumps and its upper envelope is described by the law of the iterated logarithm (\ref{lil1}); then from Theorem 8 in \cite{pa}, the reflected process  $((Y-\underline{Y})_{(T_0^- -t)^{-}}, 0\leq t\leq T_0^-)$ also satisfies the same law of the iterated logarithm (\ref{lil2}).\QED

It is important to note that when  $\alpha=2$, the upper envelope of the time reversed process 
$\{(Y_{(T_0^- -t)^{-}}, 0\leq t\leq T_0^-),\px\}$ is described by the Kolmogorov-Dvoretsky-Erd\H os integral test since the latter process has the same law as the square of a Bessel process of dimension $d=3$  killed at its last passage time below the level $x>0$. (see for instance It\^ o and McKean \cite{im})

\subsection{Some explicit calculations for $\xi$.}
Here, we compute explicitly some   functionals of the L\'evy process $\xi$ using the pssMp-Lamperti representation of $(Y, \px)$ and the results of section 3 applied to the latter.

Recall  that, under $\p_x$, $\xi$ is a spectrally positive L\'evy process starting from $0$, which drifts towards $-\infty$,  whose Laplace exponent is given by
\[
 \Psi(\theta) = m\frac{\Gamma(\theta+\alpha)}{\Gamma(\theta)\Gamma(\alpha)},
\]
and does not depend on $x$. Let us define, the last passage time of $\xi$ by
\[
D_x=\sup\big\{t\geq 0:\xi_t\geq x\big\}.
\]
\begin{proposition}Let $u<0$ and $v<u$, then
\[
P\Big(\inf_{0\leq t\leq D_u}\xi_t\geq v\Big)=\big(1-e^{v-u}\big)^{\alpha-1}.
\]
\end{proposition}
\bigskip \noindent{\it Proof:} Firstly we remark that it is known for spectrally negative stable processes of index $\alpha\in(1,2]$ that the scale function $W(x)$ is proportional to $x^{\alpha-1}$. Secondly note from Corollary \ref{noarrow}, we have  that
\[
U_y=x^{\alpha-1}\int_0^{D_{\log y/x}} e^{(\alpha-1)\xi_u} \ud u
\]
where $U_y=\sup\{t\geq 0: Y_t \geq y\} $.
Again, using the pssMp-Lamperti representation of $(Y,\px)$, we deduce
\[
P\Big(\inf_{0\leq t\leq D_u}\xi_t\geq v\Big)=\px\Big(\inf_{0\leq t\leq U_y} Y_t\leq z\Big),
\]
where $v=\log(z/x)$ and $u=\log(y/x)$ and $y\leq z$. Hence, from Corollary 1, 
we have
\[
P\Big(\inf_{0\leq t\leq D_u}\xi_t\geq v\Big)=\big(1-e^{v-u}\big)^{\alpha-1},
\]
which proves our result.\QED
\begin{proposition}  There exist a constant $k>0$ such that,
\[
\liminf_{t\to 0}\frac{\xi_t}{t^{1/\alpha}(\log\log (1/t))^{1-1/\alpha}}=-k\qquad P\textrm{-a.s.}
\]
\end{proposition}
\noindent{\it Proof:}  First, recall that $(X,\px)$ denotes a stable L\'evy process with no negative jumps starting from $x>0$. From Theorem VIII.5 in \cite{Be}, we know that there exist a constant $k>0$ such that
\[
\liminf_{t\to 0}\frac{X_t-1}{t^{1/\alpha}(\log\log (1/t))^{1-1/\alpha}}=-k\qquad \p_1\textrm{-a.s.},
\]
which implies that
\[
\liminf_{t\to 0}\frac{Y_t-1}{t^{1/\alpha}(\log\log (1/t))^{1-1/\alpha}}=-k\qquad \p_1\textrm{-a.s.},
\]
since $\p_1$-a.s., $\theta(t)\sim t$ as $t$ goes to $0$.\\
On the other hand from the pssMp-Lamperti representation of $(Y,\p_1)$, we have
\[
\liminf_{t\to 0}\frac{Y_t-1}{t^{1/\alpha}(\log\log (1/t))^{1-1/\alpha}}=
\liminf_{t\to 0}\frac{e^{\xi_{\zeta(t)}}-1}{t^{1/\alpha}(\log\log (1/t))^{1-1/\alpha}}\qquad  \p_1\textrm{-a.s.}.
\]
Next, since $\p_1$-a.s., $\zeta(t)\sim t$ as $t\to 0$, we deduce from the above identity that
\[
\liminf_{t\to 0}\frac{\xi_t}{t^{1/\alpha}(\log\log (1/t))^{1-1/\alpha}}=-k \qquad P\textrm{-a.s.},
\]
and the proof is now complete.\QED

Now, we define the first passage time of $\xi$ above $a\in \R$ by
\[
S^{+}_a=\inf\big\{t\geq 0: \xi \ge a\big\}.
\]
We also introduce the Mittag-Leffler function of parameter $\alpha$ by
\[
E_{\alpha}(x)=\sum_{n=0}^{\infty}\frac{x^n}{\Gamma(1+\alpha n)}, \qquad x\in \R,
\]
and its derivative by $E'_{\alpha}$. The following result specifies the law of some exponential functionals of the L\'evy process $\xi$.
\begin{theorem} For each $a>0$ and $q>$, we have
\begin{equation}\label{expflp1}
\begin{split}
E\left(\exp\left\{-q\int_0^{S^{+}_a}e^{\alpha\xi_s}\ud s\right\}\ind_{\{S^+_a<\infty\}}\right)&=E_{\alpha}(q(e^a-1)^\alpha)\\
&-(1-e^{-a})^{\alpha-1}E'_{\alpha}(q(e^a-1)^{\alpha})\frac{E_{\alpha}(qe^{\alpha a})}{ E'_{\alpha}(qe^{\alpha a})}, 
\end{split}
\end{equation}
and
\begin{equation}\label{expflp2}
E\left(\exp\left\{-q\int_0^{\infty}e^{\alpha\xi_s}\ud s\right\}\ind_{\{S^+_a=\infty\}}\right)=(1-e^{-a})^{\alpha-1}\frac{E'_{\alpha}(q(e^a-1)^{\alpha})}{ E'_{\alpha}(qe^{\alpha a})}.
\end{equation}
\end{theorem}
\bigskip \noindent{\it Proof:} The proof follows from Theorem 2 and the pssMp-Lamperti representation of  $(Y,\px)$. We first recall that  for $t\geq 0$  
\[
I_{t}((\alpha-1)\xi)=\int_0^{t}e^{(\alpha-1)\xi_u}\ud u,
\]
and  that $T^-_0=x^{\alpha-1}I_{\infty}((\alpha-1)\xi)$.
From the pssMp-Lamperti representation of $Y$ and the change of variabel $s=x^{\alpha-1}I_u((\alpha-1)\xi)$, we have on $\{T^+_b<T^-_0\}$ 
\[
\int_0^{T^+_b}Y_s\ud s=x\int_0^{x^{\alpha-1}I_{S^{+}_{\log b/x}}((\alpha-1)\xi)}\exp\Big\{\xi_{h(s/x^{\alpha-1})}\Big\}\ud s=x^{\alpha}\int_0^{S^{+}_{\log b/x}}\exp\Big\{\alpha\xi_{u}\Big\}\ud u.
\]
This implies that
\begin{equation}\label{tprog}
\e_x\left(\exp\left\{-q\int_0^{T^+_b}Y_s\ud s\right\}\ind_{\{T^+_b<T^-_0\}}\right)=E\left(\exp\left\{-qx^{\alpha}\int_0^{S^{+}_a}e^{\alpha\xi_s}\ud s\right\}\ind_{\{S^+_a<\infty\}}\right),
\end{equation}
with $a=\log(b/x)$.

On the other hand, from the proof of Theorem 1 in Bertoin \cite{Be1}, it is known 
\[
W^{(q)}(y)=\alpha y^{\alpha-1}E'_{\alpha}(qy^{\alpha})\qquad\textrm{ and }\qquad Z^{(q)}(y)=E_{\alpha}(qy^{\alpha}) \qquad \textrm{for}\quad y\geq 0.
\] 
Now, from (\ref{tprog}), the above formulas and applying Theorem 2, we deduce (\ref{expflp1}).

The equality (\ref{expflp2}) follows by similar arguments, we just need to note that
\[
\int_0^{T_0^-}Y_s\ud s=x^{\alpha}\int_0^{\infty}\exp\Big\{\alpha\xi_{u}\Big\}\ud u.
\]
The proof is now complete.\QED
\section{Self-similar CBI-processes.}
In light of the previous section, we may also consider the process $(Y, \px^\uparrow)$ in the context of self-similarity. To this end, let us recall that $\xi^*=(\xi^*_t : t\geq 0)$ is a L\'evy process, under $\px$, which drift towards $\infty$ and whose Laplace exponent is given by
\[
\Psi^*(\theta) = m^* \frac{\Gamma(\theta + \alpha - 1)}{\Gamma(\theta - 1)\Gamma(\alpha)}, \qquad\textrm{for }\quad\theta \geq 0
\]
where $m^*=\mathbb{E}(\xi_1^*)$. Also recall that when $\alpha=2$, we have that $\xi^*_t=B_t+t/2$, where $B$ is a standard Brownian motion.
\begin{corollary}
The process $(Y,\px^\uparrow)$ is a positive self-similar Markov process with index of self-similarity $(\alpha-1)$. Moreover, its pssMp-Lamperti representation under $\mathbb{P}^{\uparrow}_x$ is given by
\[
Y_t = x\exp\{\xi^*_{\zeta^*(tx^{-(\alpha-1)})}\}.
\]
where 
\[
\zeta^*(t) = \inf\left\{ s\geq 0 : \int_0^s \exp{\Big\{(\alpha -1)\xi^*_u \Big\}}\ud u > t \right\}.
\]
Moreover, the process $(Y,\p^{\uparrow}_x)$ converges weakly as $x$ tend to $0$, in the sense of Skorokhod towards $(Y,\p^{\uparrow})$, a pssMp starting from $0$ with same semigroupe as $(Y,\px^{\uparrow})$, for $x>0$, and with entrance law given by
\begin{equation}\label{entlci}
\e^{\uparrow}\Big(e^{-\lambda Y_t}\Big)=\Big(1+c_+(\alpha-1)t\lambda^{\alpha-1}\Big)^{-\frac{\alpha}{(\alpha-1)}}.
\end{equation}
\end{corollary}
It is important to note that the process $(Y,\p^{\uparrow})$ can be constructed as in Theorem 1 in Caballero and Chaumont \cite{cc1}.

\bigskip \noindent{\it Proof:} The pssMp-Lamperti representation follows from Proposition \ref{generic} where now the underlying L\'evy process is $\xi^*$. Now, note that the L\'evy process $\xi^*$ satisfies the conditions of Theorem 1 and 2 in \cite{cc1}. Hence the family of processes $\big\{(Y_t, t\geq 0), \px^{\uparrow}\big\}$, for $x>0$, converges weakly with respect to the Skorohod topology, as $x$ goes to $0$, towards a pssMp starting from 0 which is $(Y,\p^{\uparrow})$. It is well-known (see for instance \cite{lamb0}) that its entrance law is of the form 
\[
\e^{\uparrow}\Big(e^{-\lambda Y_t}\Big)=\exp\left\{-\int_0^{t}\phi(u_s(\lambda))\ud s\right\}.
\]
Solving (\ref{DEut}) explicitly we find that 
\[
u_t(\lambda) = [c_+(\alpha-1)t + \lambda^{-(\alpha-1)}]^{-\frac{1}{(\alpha-1)}}.
\]
We obtain (\ref{entlci}) from straightforward calculations, recalling that $\phi(\lambda) = c_+\alpha\lambda^{\alpha-1}$ thus completing the proof. \QED 
This last  theorem, Theorem 2 and proposition 1 give rise to the following flow of transformations. Let $h_1(y)=e^y, h_2(y)=e^{-y}, h_3(y)=y$ and $h_{4}=1/y$, then
\[
\hspace{-0.5cm}
\begin{array}{ccccccc}
\begin{array}{c} \xi \\ \end{array}& {\begin{array}{c} \text{\footnotesize pssMp-Lamp} \\  {\leftarrow -- \rightarrow} \end{array}}&\begin{array}{c} (X,\px) \\ \end{array} & {\begin{array}{c} \text{\footnotesize CB-Lamp} \\  {\leftarrow -- \rightarrow} \end{array}} &\begin{array}{c} (Y,\px) \\ \end{array} &{\begin{array}{c} \text{\footnotesize pssMp-Lamp} \\  {\leftarrow -- \rightarrow} \end{array}}  &\begin{array}{c} \xi \\ \end{array}\\
{\begin{array}{ccc} & \uparrow & ^{h_2} \\& | & \\ & |&  \\ & \downarrow &_{h_1}\end{array}} & \alpha&
{\begin{array}{ccc} ^{h_4}& \uparrow &  \\& | & \\ & |&  \\ _{h_3}& \downarrow & \end{array}} & & {\begin{array}{ccc} & \uparrow & ^{h_4} \\& | & \\ & |&  \\ & \downarrow & _{h_3}\end{array}} &\alpha-1 &{\begin{array}{ccc} ^{h_2}& \uparrow &  \\& | & \\ & |&  \\ _{h_1}& \downarrow & \end{array}}  \\
\begin{array}{c} \\ \xi^*  \end{array}&{\begin{array}{c} \text{\footnotesize pssMp-Lamp} \\  {\leftarrow -- \rightarrow} \end{array}} & \begin{array}{c} \\(X,\px^\uparrow)  \end{array} & {\begin{array}{c} \text{\footnotesize CB-Lamp} \\  {\leftarrow -- \rightarrow} \end{array}}  & \begin{array}{c} \\(Y,\px^\uparrow)  \end{array} &{\begin{array}{c} \text{\footnotesize pssMp-Lamp} \\  {\leftarrow -- \rightarrow} \end{array}} & \begin{array}{c} \\ \xi^*  \end{array}
\end{array}
\]
where the vertical arrows are the result of a Doob $h$-transform with the $h$-function indicated in each dircection and the parameters $\alpha$ and $\alpha-1$ are the index of self-similarity on the pssMp-Lamperti representation.

From the form of the entrance law of $(Y,\p^{\uparrow})$ and the Cram\'er's condition of the L\'evy process $\xi^*$ (see Proposition 1), we deduce the following corollary.
\begin{corollary}\label{corr3}
The exponential functional of $\xi^*$, defined by
\[
I^*:=\int_0^{\infty}e^{-(\alpha-1)\xi^*_s}\ud s,
\]
satisfies
\begin{equation}\label{expn}
E\Big(e^{-\lambda (I^*)^{-1}}(I^*)^{-1}\Big)=m^*(\alpha-1)\big(1+c_+(\alpha-1)\lambda^{\alpha-1}\big)^{-\alpha/(\alpha-1)},\quad\textrm{for } \quad\lambda\geq 0.
\end{equation}
Moreover, we have that
\[
P\big(I^*> t\big)\sim C t^{-1/(\alpha-1)} \qquad \textrm{as }\quad t\to\infty 
\]
where $C$ is a nonnegative constant which depends on $\alpha$.
\end{corollary} 
\bigskip \noindent{\it Proof:} The equality   (\ref{expn}) follows from the form of the entrance law of positive self-similar Markov process (see for instance Proposition 3 in \cite{cc1}), which says in our particular case that
\[
\e^{\uparrow}\Big(e^{-\lambda Y_1}\Big)=\frac{1}{m^*(\alpha-1)}E\Big(e^{-\lambda (I^*)^{-1}}
(I^*)^{-1}\Big),\quad \textrm{for }\quad \lambda\geq 0.
\]
According to Lemma 4 in Rivero \cite{ri},  if $\xi$ is a non arithmetic L\'evy process which drift towards $-\infty$ and satisfying the Cram\'er's condition for some $\theta>0$, i.e. $E(\exp{\theta\xi_1})=1$ implies that for $\beta>0$
\[
P\left(\int_0^{\infty}e^{\beta \xi_s}\ud s>t\right)\sim Kt^{-\theta/\beta}\qquad \textrm{as }\quad t\to\infty,
\] 
where $K$ is a nonnegative constant which depends on $\beta$. {

Hence, recalling Proposition \ref{prop1}, one may  apply Lemma 4 in \cite{ri} for the process $-\xi^*$ and get
\[
P\big(I^*> t\big)\sim C t^{-1/(\alpha-1)} \qquad \textrm{as }\quad t\to\infty, 
\]
where $C$ is a nonnegative constant which depends on $\alpha$.\QED 
The estimation of the tail behaviour of $I^*$ allow us to describe the lower envelope at $0$ and at $+\infty$ of $(Y,\p^{\uparrow})$.
\begin{theorem}
Let $f$ be an increasing function such that $\lim_{n\to 0}f(t)/t=0$, then:
\[
\p^{\uparrow}\big(Y_t<f(t), \textrm{ i.o., as } t\to 0\big)=0\textrm{ or } 1,
\]
accordingly as 
\[
\int_{0+} f(t) t^{\alpha/(\alpha-1)}\ud t \quad \textrm{ is finite or infinite}.
\]
In particular, 
\[
\liminf_{t\to 0}\frac{Y_{t}}{t^{\kappa}}=
\left\{ \begin{array}{ll}
0 &\textrm{if $\kappa>\frac{1}{\alpha-1}$}\\
+\infty &\textrm {if $\kappa\leq\frac{1}{\alpha-1}$}
\end{array} \right.
\qquad \p^{\uparrow}-\textrm{a.s.}
 \]
Let $g$ be an increasing function such that $\lim_{n\to \infty}g(t)/t=0$, then for all $x\geq 0$:
\[
\p^{\uparrow}_x\big(Y_t<g(t), \textrm{ i.o., as } t\to \infty\big)=0\textrm{ or } 1,
\]
accordingly as 
\[
\int^{\infty} \left(\frac{g(t)}{t}\right)^{1/(\alpha-1)}\frac{\ud t}{t} \quad \textrm{ is finite or infinite}.
\]
In particular, for any $x\geq 0$ 
\[
\liminf_{t\to \infty}\frac{Y_{t}}{t^{\kappa}}=
\left\{ \begin{array}{ll}
0 &\textrm{if $\kappa<\frac{1}{\alpha-1}$}\\
+\infty &\textrm {if $\kappa\geq \frac{1}{\alpha-1}$}
\end{array} \right.
\qquad \p_x-\textrm{a.s.}
 \]
\end{theorem}
\bigskip \noindent{\it Proof:} This result follows from the tail behaviour of $I^*$ described in Corollary \ref{corr3} and Theorem 3 in \cite{cp}.\QED
Now, we define the  family of positive self-similar Markov process $\widehat{X}^{(x)}$ whose pssMp-Lamperti representation is given by 
\[
X^{(x)}=\Big(x\exp\Big\{\widehat{\xi}^*_{\widehat{\zeta}^*(t/x^{\alpha-1})}\Big\}, 0\leq t\leq x^{\alpha-1}I(\xi^*)\Big), \qquad x>0,
\]
where $\widehat{\xi^*}=-\xi^*$ and 
\[
\widehat{\zeta}^*(t)=\inf\left\{t:\int_0^s\exp\Big\{(\alpha-1)\widehat{\xi^*}_u\Big\}\ud u >t \right\}.
\]
We emphasize that the random variable $ x^{\alpha-1}I^*$ corresponds to the first time at which the process $\widehat{X}^{(x)}$ hits $0$, moreover for each $x>0$, the process $\widehat{X}^{(x)}$ hits $0$ continuously.

We now set
\[
U^{-}_y=\sup\{t\geq 0: Y_t\leq y\} \qquad\textrm{ and }\qquad\Gamma=Y_{U^{-}_y}.
\]
According to Proposition 1 in \cite{cp}, the law of the process $\widehat{X}^{(x)}$ is a regular version of the law of the process $\{(Y_{(U^-_y-t)^-}, 0\leq t \leq U^{-}_y),  \p^{\uparrow}\}$
conditionally on $\{\Gamma=x\}$, $x\in [0,y]$. Hence, the latter process  is equal in law to 
\[
\Big(\gamma \exp\Big\{\widehat{\xi}^{*}_{\widehat{\zeta}^{*}(t/\Gamma^{\alpha-1})}\Big\}, 0\leq t\leq \gamma^{\alpha-1}I^*\Big)
\]
where $\gamma$ is equal in distribution to $\Gamma$ and independent of $\widehat{\xi}^*$.

We shall momentarily turn out attention to describing the law of $\Gamma$. Let  $H=(H_t,t\geq 0)$ be the ascending ladder height process associated to $\xi^*$ (see Chapter VI in \cite{Be} for a formal definition) and denote by $\nu$ its L\'evy measure. Hence according to Lemma 1 in \cite{cp}, the law of $\Gamma$ is characterized as follows:
\[
\log(y^{-1}\Gamma)\ed-\mathcal{U}\mathcal{Z},
\]
where $\mathcal{U}$ and $\mathcal{Z}$ are independent r.v.'s, $\mathcal{U}$ is uniformly distributed over $[0,1]$ and the law of $\mathcal{Z}$ is given by:
\[
P(\mathcal{Z}>u)=E(H_1)^{-1}\int_{(u,\infty)}s\nu(\ud s), \qquad u\geq 0.
\]
The above discussion now gives us the following propositions.
\begin{proposition} Let $z\geq y>0$, then
\[
\p^{\uparrow}\Big(\sup_{0\leq s\leq \sigma_y} X_s\leq z\Big)=\p^{\uparrow}\Big(\sup_{0\leq s\leq U^-_y} Y_s\leq z\Big)=1-\frac{y}{m^*z}.
\]

\end{proposition}
\bigskip \noindent{\it Proof:} The first equality follows from the CB-Lamperti representation of $(Y,\p^{\uparrow})$. Now, from the time reversal of $(Y,\p^{\uparrow})$ at its last passage time $U^{-}_y$, we have that
\[
\p^{\uparrow}\Big(\sup_{0\leq s\leq U^-_y} Y_s\leq z\Big)=P\Big(\sup_{s\geq 0 } \widehat{\xi}^*_s\leq \log(z/\gamma)\Big)
\]
where $\gamma$ is equal in distribution to $\Gamma$ and independent of $\widehat{\xi}^*$.
On the one hand, it is a well established fact that the all-time supremum of a spectrally negative L\'evy process which drifts to $-\infty$ is exponentially distributed with parameter equal to the largest root of its Laplace exponent. In particular,  for $x\geq 0$,
\[
P\Big(\sup_{s\geq 0 } \widehat{\xi}^*_s\leq x\Big)=1-e^{-x}.
\]
Note that by inspection of $\Psi^*(\theta)$ the largest root is clearly $\theta=1$ (there are at most two and one of them is always $\theta=0$).
On the other hand, from the above discussion, the random variables $\widehat{\xi^*}$ and $\Gamma$ are independent. Hence
\[
\p^{\uparrow}\Big(\sup_{0\leq s\leq U^-_y} Y_s\leq z\Big)=E\left(1-\frac{\Gamma}{z}\right).
\]
Therefore, in order to complete the proof, it is enough to show that $E(\Gamma)=\displaystyle\frac{y}{m^*}.$

To this end, recall that $\Gamma\ed ye^{\mathcal{U}\mathcal{Z}}$, where $\mathcal{U}$ and $\mathcal{Z}$ are independent r.v.'s and defined above,  hence
\begin{eqnarray}
E(\Gamma)&=&y\int_{(0,\infty)}\int_0^1 e^{-uz}\ud u P(\mathcal{Z}\in \ud z)\notag \\
&=&y\int_{(0,\infty)}\frac{1}{z}(1-e^{-z}) P(\mathcal{Z}\in \ud z)\notag\\
&=&\frac{y}{E(H_1)}\int_{(0,\infty)}(1-e^{-z})\nu(\ud z).\label{pieces}
\end{eqnarray}
Next note that since $\Psi^*(\theta)$ has its largest root at $\theta=1$, the Wiener-Hopf factorization for the process $\widehat{\xi}^*$ must necessarily take the from $\Psi^*(\theta) =(\theta-1)\phi(\theta)$ for $\theta\geq 0$, where $\phi(\theta)$ is the Laplace exponent of the descending ladder height process of $\widehat{\xi}^*$. Note that $\phi$ has no killing term (i.e. $\phi(0)=0$) as $\widehat{\xi}^*$ drifts to $-\infty$. Moreover, $\phi$ has no drift term as $\widehat{\xi}^*$ has no Gaussian component (cf. \cite{ckp}). The latter two observations imply that 
\[
\int_{(0,\infty)} (1- e^{-z})\nu(\ud z)=\phi(1) = \left.\frac{\Gamma(\theta-1+\alpha)}{(\theta-1)\Gamma(\theta-1)\Gamma(\alpha)}\right|_{\theta=1} =1.
\]
Note also that $-m^* = \mathbb{E}(\widehat{\xi}^*_1)=\Psi^{*\prime}(0+) = - \phi'(0+)= \mathbb{E}(H_1)$. Putting the pieces together in (\ref{pieces}) completes the proof.
\QED

Now, we define the following exponential functional of $\xi^*$,
\[
I':=\int_{0}^{\infty} e^{-\alpha \xi^*_s}\ud s.
\]
The exponential functional $I'$ was studied by Chaumont et al. \cite{ckp}. In particular the authors in \cite{ckp} found that \[
P\big(1/I'\in \ud y\big)=\alpha m^*q_1(y)\ud y,
\]
where $q_1$ is the density of the entrance law of the measure of the excursions away from $0$ of the reflected process $(X_t-\inf_{s\leq t}X_s, t\geq 0)$, under $\p$.

The time reversal property of $(Y,\p^{\uparrow})$ at its last passage time combined with the CB-Lamperti representation and the pssMp-Lamperti representation give us the following result for the total progeny of the self-similar CB-process with immigration.
\begin{proposition}
The total progeny of $(Y,\p^{\uparrow})$, the self-similar CB-process with immigration starting from $0$, up to time $U^-_y$, for $y>0$, i.e.
\[
J_{U^-_y}:=\int_0^{U^-_y}Y_s\ud s , \qquad \textrm{under }\quad \p^{\uparrow},
\] 
and $\sigma_y$ the last passage time of  $(X,\p^{\uparrow})$, the stable process conditioned to stay positive, below $y$ are both equal in law to
$\Gamma^{\alpha}I'$.

\end{proposition}

\section{Quasi-stationarity}

We conclude this paper with some brief remarks on a different kind of conditioning of CB-processes to (\ref{CBCSP}) which results in a so-called quasi-stationary distribution for the special case of the self-similar CB-process.
Specificially we are interested in establishing the existence of normalization constants $\{c_t : t\geq 0\}$ such that the weak limit
\[
 \lim_{t\uparrow\infty}\mathbb{P}_x(Y_t/c_t\in {\rm d}z|T^-_0>t)
\]
exists for $x> 0$ and $z\geq 0$.
 
Results of this kind have been established for CB-processes for which the underlying spectrally positive L\'evy process has a Gaussian process in \cite{lamb} and \cite{li}. In the more general setting, \cite{pakes} formulates conditions for the existence of such a limit and characterizes the resulting quasi-stationary distribution. The result below shows that in the self-similar case we consider in this paper, an explicit formulation of the normalization sequence $\{c_t: t\geq 0\}$ and the limiting distribution is possible. 

\begin{lemma}Fix $\alpha\in(1,2]$.
 For all $x \geq 0$, with $c_t= [c_+(\alpha-1)t]^{1/(\alpha-1)}$ 
\[
 \lim_{t\uparrow\infty}\mathbb{E}_x(e^{-\lambda Y_t/c_t}|T^-_0>t) = 1 - \frac{1}{[1+ \lambda^{-(\alpha-1)}]^{1/(\alpha-1)}}.
\]
\end{lemma}

\bigskip \noindent{\it Proof:}
The proof pursues a similar line of reasoning to the the aforementioned references \cite{lamb, li, pakes}.
From (\ref{CBut}) it is straightforward to deduce that 
\[
 \lim_{t\uparrow\infty}\mathbb{E}_x(1- e^{-\lambda Y_t/c_t}|T^-_0>t) = \lim_{t\uparrow\infty}\frac{u_t(\lambda/c_t)}{u_t(\infty)}
\]
if the limit on the right hand side exists. However, since $\psi(\lambda) = \lambda^\alpha$
it is easily deduced from (\ref{DEut}) that 
\[
 u_t(\lambda) = [c_+(\alpha-1)t + \lambda^{-(\alpha-1)}]^{-1/(\alpha-1)}
\]
and the result follows after a brief but easy calcluation.
\QED

Although quasi-sationarity in the sense of `conditioning to stay positive' does not make sense in the case of the CBI-process $(Y,\mathbb{P}_x^\uparrow)$, it appears that the normalizing constants $\{c_t: t\geq 0\}$ serve a purpose to obtain the below distributional limit. A similar result is obtained in \cite{lamb} for CBI-processes with a Gaussian component in the underlying spectrally positive L\'evy process.

\begin{lemma}
 Fix $\alpha\in(1,2]$.
 For all $x \geq 0$, with $c_t= [c_+(\alpha-1)t]^{1/(\alpha-1)}$ 
\[
 \lim_{t\uparrow\infty}\mathbb{E}^\uparrow_x(e^{-\lambda Y_t/c_t}) = \frac{1}{[ \lambda^{(\alpha-1)} +1]^{\alpha/(\alpha-1)}}.
\]
\end{lemma}
\bigskip \noindent{\it Proof:} We follow ideas found in Lambert \cite{lamb}.
In the latter paper, it is shown that $(Y, \mathbb{P}^\uparrow_x)$ may also be obtained as the Doob $h$-transform of the process $(Y, \mathbb{P})$ with $h(x)=x$. That is to say
\[
 \mathbb{E}^\uparrow_x( e^{-\lambda Y_t})=  \mathbb{E}^\uparrow_x(Y_t e^{-\lambda Y_t}).
\]
Differentiating (\ref{CBut}) this implies that 
\[
  \mathbb{E}^\uparrow_x( e^{-\lambda Y_t}) = e^{-x u_t(\lambda)}\frac{\psi(u_t(\lambda))}{\psi(\lambda)}.
\]
Plugging in the necessary expressions for $\psi$ and $u_t(\lambda)$ as well as replacing $\lambda$ by $\lambda/c_t$ in the previous formula, the result follows directly.
\QED

Note the limiting distribution has is the law at unit time of a Gamma subordinator with parameters 1 and $\alpha/(\alpha-1)$ time changed by a stable subordinator of index $\alpha-1$. This is also equal to the law of $Y_1$ under $\mathbb{P}^\uparrow$ when $c_+ = 1/(\alpha-1)$.

\section*{Acknowledgments} Both authors would like to thank Victor Rivero for many useful discussions.
This research was funded by EPSRC grant EP/D045460/1.

\end{document}